\documentclass[11pt, reqno]{amsart}
\usepackage[utf8]{inputenc}
\usepackage{amsmath}
\usepackage{amsthm}
\usepackage{amsfonts}
\usepackage{amssymb}
\usepackage{mathrsfs, mathtools, mathdots}
\usepackage{diagbox}

\usepackage[top=27mm,bottom=27mm,left=27mm,right=27mm]{geometry}
\newtheorem{theorem}{Theorem}[section]

\newtheorem{corollary}[theorem]{Corollary}

\newtheorem*{theorem*}{Theorem}
\newtheorem*{remark*}{Remark}
\newtheorem*{problem*}{Problem}
\newtheorem*{lemma*}{Lemma}
\newtheorem*{conjecture*}{Conjecture}
\newtheorem{lemma}[theorem]{Lemma}
\newtheorem{prop}[theorem]{Proposition}

\begin{document}

\title[Transcendence of special values of modular functions]{A note on transcendence of special values \\ of functions related to modularity}

\author[T. Bhowmik]{Tapas Bhowmik}

\author[S. Pathak]{Siddhi Pathak}

\address{Department of Mathematics, University of South Carolina, 1523 Greene Street, LeConte College, Columbia, South Carolina, USA 29208.}

\address{Chennai Mathematical Institute, H-1 SIPCOT IT Park, Siruseri, Kelambakkam, Tamil Nadu, India 603103.}

\email{tbhowmik@email.sc.edu}

\email{siddhi@cmi.ac.in}



\thanks{Research of the second author was partially supported by an INSPIRE Faculty fellowship.}

\begin{abstract}
In this note, we study the arithmetic nature of values of modular functions, meromorphic modular forms and meromorphic quasi-modular forms with respect to arbitrary congruence subgroups, that have algebraic Fourier coefficients. This approach unifies many of the known results, and leads to generalizations of the theorems of Schneider, Nesterenko and others. 
\end{abstract}

\maketitle

\section{\bf Introduction}
\bigskip

The study of the arithmetic nature of values of special transcendental functions at algebraic arguments has been a well-established theme in number theory. Continuing in this spirit, this note focuses on the transcendental nature and algebraic independence of values of functions arising in the modular world, such as modular functions, modular forms and quasi-modular forms (see Section 2 for definition of the functions appearing below). The genesis of this study can be traced back to a 1937 theorem of Schneider \cite{schneider}, namely,
\begin{theorem}[Schneider]\label{schneider}
    If $\tau\in \mathbb{H}$ is algebraic but not imaginary quadratic, then $j(\tau)$ is transcendental.
\end{theorem}
It is known from the theory of complex multiplication that if $\tau\in\mathbb{H}$ generates an imaginary quadratic field ($\tau$ is a CM point), then $j(\tau)$ is an algebraic number. Therefore, Schneider's theorem translates to the statement: if $\tau$ is algebraic, then $j(\tau)$ is algebraic if and only if $\tau$ is CM. A further conjecture by Mahler \cite{mahler} in this regard, proved by Barr\'e-Sirieix, Diaz, Gramain and Philibert \cite{b-d-g-p} states that
\begin{theorem}[Barr\'e-Sirieix, Diaz, Gramain, Philibert]\label{mahl}
For any $\tau\in \mathbb{H}$, at least one of the two numbers $e^{2\pi i\tau}$ and $j(\tau)$ is transcendental.
\end{theorem}
This can be derived as a consequence of a remarkable theorem of Nesterenko \cite{nesterenko} with several applications.
\begin{theorem}[Nesterenko] \label{nest}
If $\tau\in \mathbb{H}$, then at least three of the numbers
$$
e^{2 \pi i \tau},\quad E_2(\tau),\quad E_4(\tau),\quad E_6(\tau)
$$
are algebraically independent over $\overline{\mathbb{Q}}$.
\end{theorem}

Although the above theorems are about specific functions of ``level $1$", it is the aim of this note to highlight that they are sufficient to deduce the corresponding results for functions associated with arbitrary congruence subgroups. The authors believe that this fact may be known to experts, but is not well-documented. Often, the congruence subgroup in question is restricted to be the group $\Gamma_0(N)$. The results in this paper apply to functions satisfying appropriate modularity properties with respect to $\Gamma(N)$, and hence, arbitrary congruence subgroups.\\

A detailed investigation of algebraic independence of values of modular forms and quasi-modular forms was carried out by S. Gun, M. R. Murty and P. Rath \cite{r-m-g} in 2011. Their results on values of modular forms were further elaborated upon for higher level in \cite{hamieh-murty} by A. Hamieh and M. R. Murty\footnote{A small correction in their statement of Theorem 1.2 is required. The conclusion should read as ${(\pi/\omega_{\tau})}^k \, L_q(k,\chi)$ is algebraic, without the $L(1-k,\chi)$ term.}. These theorems will follow from our discussion later. In the context of quasimodular forms, it was proven independently by Gun, Murty and Rath \cite[Theorem 7]{r-m-g} and C. Y. Chang \cite{chang} that
\begin{theorem}[Gun-Murty-Rath and Chang]\label{GMR-Chang}
    If $f$ is a quasi-modular form of non-zero weight for $SL_2(\mathbb{Z})$ with algebraic Fourier coefficients and $\tau\in\mathbb{H}$ is such that $j(\tau)\in\overline{\mathbb{Q}}$, then either $f(\tau) = 0$ or $f(\tau)$ is transcendental.  
\end{theorem}
The above statement is also proved for quasimodular forms with respect to $\Gamma_0(N)$ in \cite{r-m-g}.\\

Another instance of investigation is a recent paper of D. Jeon, S.-Y. Kang and C. H. Kim \cite[Theorem 2.4]{jeon-kang-kim}, where they prove the following. Let $N \in \mathbb{N}$, $\mathfrak{g}:=\mathfrak{g}_0(N) = $ genus of $X_0(N)$, the modular curve obtained as the quotient of the extended upper half plane by $\Gamma_0(N)$. Suppose that $\mathfrak{g} > 0$ and $m \geq \mathfrak{g}+1 $ is an integer. Let $\mathfrak{f}_{N,m}$ denote the unique modular function with respect to $\Gamma_0(N)$ constructed in \cite{j-k-k-2} such that $\mathfrak{f}_{N,m}$ is holomorphic on the upper half plane and
\begin{equation*}
    \mathfrak{f}_{N,m}(q) = \frac{1}{q^m} + \sum_{l=1}^\mathfrak{g} a_N(m,-l) \frac{1}{q^l} + O(q),
\end{equation*}
with the coefficients of powers of $q$ being algebraic numbers. Then, they show the following.
\begin{theorem}[Jeon, Kang, Kim]\label{jjk-thm}
   Let $f$ be a non-zero meromorphic modular form with respect to $\Gamma_0(N)$ with algebraic Fourier coefficients. If $\tau$ is either a zero or a pole of $f$, then $\mathfrak{f}_{N,m}(\tau)$ is algebraic for all $m$. 
\end{theorem}
As a corollary, they deduce that any zero or pole of $f$ should be either CM or transcendental.\\

In this note, we first give an exposition of the algebraic structure of modular functions of higher level, following Shimura \cite{shimura}. Building upon this and using Schneider's theorem, we prove
\begin{theorem}\label{mdf-values-thm}
Let $g$ be a non-constant modular function with respect to a congruence subgroup $\Gamma$ of level $N$, with algebraic Fourier coefficients at $i\infty$. 
\begin{enumerate}
    \item[(a)] If $\tau\in\mathbb{H}$ is either a zero or a pole of $g$, then $j(\tau)$ is algebraic, and hence, $\tau$ is either CM or transcendental.
    \item[(b)] If $\tau\in\mathbb{H}$ is such that $\tau$ is not a pole of $g$, then $j(\tau)\in\overline{\mathbb{Q}}\Leftrightarrow g(\tau)\in\overline{\mathbb{Q}}$. Thus, at least one of $g(\tau)$ and $e^{2\pi i\tau}$ is transcendental for any $\tau \in \mathbb{H}$.
\end{enumerate}
\end{theorem}

In the context of meromorphic modular forms, we establish a generalization of Theorem \ref{jjk-thm} to arbitrary modular forms and arbitrary modular functions. 
\begin{theorem}\label{zero-pole-value}
Let $f$ be a non-constant meromorphic modular form of weight $k \in \mathbb{Z}$ with respect to a congruence subgroup of level $N$ and $g$ be a non-constant modular function with respect to a congruence subgroup of level $M$. Suppose that both $f$ and $g$ have algebraic Fourier coefficients at $i\infty$. Let $\tau \in \mathbb{H}$ be a zero or a pole of $f$. Then either $g(z)$ has a pole at $z=\tau$ or $g(\tau)$ is algebraic.
\end{theorem}

Furthermore, we generalize the theorems in \cite{r-m-g} and \cite{hamieh-murty} to the setting of meromorphic modular forms. In particular, we show the following. 
\begin{theorem}\label{general-mero-values-thm}
Let $f$ be a non-constant meromorphic modular form with respect to a congruence subgroup $\Gamma$ of level $N$. Suppose that $f$ has algebraic Fourier coefficients at $i\infty$. Suppose that $\tau$ is not a pole of $f$. 
  \begin{enumerate}
      \item[(a)] If $\tau\in\mathbb{H}$ is such that $e^{2\pi i\tau}$ is algebraic, then $f(\tau)$ is transcendental.
      \item[(b)] If $\tau\in\mathbb{H}$ is such that $j(\tau)\in \overline{\mathbb{Q}}$, then there exists a transcendental number $\omega_\tau$ which depends only on $\tau$ and is $\overline{\mathbb{Q}}$-linearly independent with $\pi$, such that $\left(\dfrac{\pi}{\omega_\tau}\right)^k f(\tau)\in \overline{\mathbb{Q}}$. Therefore, $f(\tau)$ is either zero or transcendental.
  \end{enumerate}  
\end{theorem}

With regard to quasi-modular forms, we extend the previously known results to meromorphic quasi-modular forms and prove
\begin{theorem}\label{quasi-mdf-thm}
    Let $\widetilde{f}$ be a non-constant meromorphic quasi-modular form with depth $p \geq 1$, with respect to a congruence subgroup $\Gamma$ of level $N$. Suppose that $\widetilde{f}$ has algebraic Fourier coefficients at $i\infty$ and that $\tau$ is not a pole of $\widetilde{f}$. 
  \begin{enumerate}
      \item[(a)] If $\tau\in\mathbb{H}$ is such that $e^{2\pi i\tau}$ is algebraic, then $\widetilde{f}(\tau)$ is transcendental.
      \item[(b)] Let 
      \begin{equation*}
          \widetilde{f} = \sum_{r=0}^p f_{r} \, E_2^r \qquad \text{ with } \qquad f_{r} \in M^m_{k-2r, N, \overline{\mathbb{Q}}}.
      \end{equation*}
      Here $M^m_{j,N,\overline{\mathbb{Q}}}$ denotes the space of meromorphic modular forms of weight $j$, level $N$ and algebraic Fourier coefficients at $i\infty$. If $\tau\in\mathbb{H}$ is such that $j(\tau)\in \overline{\mathbb{Q}}$, then 
      \begin{equation*}
          \widetilde{f}(\tau) = 0 \qquad \iff \qquad f_r(\tau) = 0 \text{  for all  } 0 \leq r \leq p.
      \end{equation*}
      Moreover, if $\widetilde{f}(\tau) \neq 0$, then $\widetilde{f}(\tau)$ is transcendental.
  \end{enumerate}  
\end{theorem}
This generalizes Theorem \ref{GMR-Chang} as well as \cite[Theorem 1.8]{balu-gun}.  \\

Let $\mathcal{Z}_{j, \overline{\mathbb{Q}}}:= \left\{ \tau \in \mathbb{H} \,: \, j(\tau) \in \overline{\mathbb{Q}} \right\}$. From the above results, we deduce the following interesting corollary.
\begin{corollary}
Let $\mathcal{Z}_{\text{mdfn}}$ be the set of zeros and poles of modular \emph{functions} of arbitrary level with algebraic Fourier coefficients, $\mathcal{Z}_{\text{mdfrm}}$ be the set of zeros and poles of meromorphic modular forms of arbitrary level with algebraic Fourier coefficients and $\mathcal{Z}_{\text{quasi-mdf}}$ be the set of zeros and poles of meromorphic quasi-modular forms of arbitrary level with algebraic Fourier coefficients. Then
\begin{equation*}
    \mathcal{Z}_{\text{quasi-mdf}} \subseteq \mathcal{Z}_{\text{mdfrm}} \subseteq \mathcal{Z}_{\text{mdfn}} \subseteq \mathcal{Z}_{j, \overline{\mathbb{Q}}}.
\end{equation*}
In particular, zeros and poles of quasi-modular forms, modular forms and modular functions are either CM or transcendental.
\end{corollary}

Finally, we have the following generalization of Theorem \ref{nest}. 
\begin{theorem}\label{tr-deg-thm}
Let $f$ be a non-constant meromorphic modular form of weight $k \in \mathbb{Z}$ with respect to a congruence subgroup of level $N$, $g$ be a non-constant modular function with respect to a congruence subgroup of level $M$ and $\widetilde{f}$ be a non-constant meromorphic quasi-modular function of depth at least $1$, with respect to a subgroup of level $\widetilde{N}$. Suppose that $f$, $g$ and $\widetilde{f}$ have algebraic Fourier coefficients at $i\infty$. If $\tau\in\mathbb{H}$ is such that $f(\tau)\neq 0$, $\widetilde{f}(\tau)\neq 0$ and $\tau$ is not a pole of $f$, $g$ and $\widetilde{f}$, then 
    $$
\operatorname{trdeg}_{\,\mathbb{Q}}\,  \,\overline{\mathbb{Q}}\left( e^{2\pi i\tau},\, f(\tau),\, g(\tau),\, \widetilde f(\tau)\right)\geq 3. 
    $$
\end{theorem}
This theorem is in the same spirit as \cite[Theorem 1.2]{wang}, where W. Wang considers the algebraic independence of the values of three algebraically independent quasi-modular forms. However, Theorem \ref{tr-deg-thm} allows one to also compare values of modular functions with those of quasi-modular forms. It can be shown that a modular function, a modular form of non-zero weight and a quasi-modular form of positive depth are algebraically independent. We include a proof of this assertion in Theorem \ref{algebraic-indep-fn-thm} for completeness. \\

\noindent We also remark that in Theorem \ref{zero-pole-value} and Theorem \ref{tr-deg-thm}, one can replace a meromorphic modular form by a half-integer weight modular form with algebraic Fourier coefficients since the square of a half-integer weight modular form is an integer weight modular form.

\bigskip

\section{\bf Preliminaries}
\bigskip

The aim of this section is to study the algebraic structure of the field of modular functions and to record the required results from transcendental number theory. For the sake of completeness and clarity of exposition, a brief account of the proofs is included, and appropriate references are given.

\subsection{Modular and quasi-modular forms}
For each $N\in \mathbb{N}$, $$
\Gamma(N):=\left\{  \left( 
\begin{array}{cc}
  a & b \\
  c & d
\end{array}
\right)\in SL_2(\mathbb{Z}):\left( 
\begin{array}{cc}
  a & b \\
  c & d
\end{array}
\right)\equiv \left( 
\begin{array}{cc}
  1 & 0 \\
  0 & 1
\end{array}
\right) (\operatorname{mod}N) \right\},
$$ 
with $\Gamma(1)=SL_2(\mathbb{Z})$. A subgroup $\Gamma$ of $SL_2(\mathbb{Z})$ is said to be congruence subgroup if there exists $N\in \mathbb{N}$ such that $\Gamma(N)\subseteq \Gamma$. The smallest such $N$ is said to be the level of $\Gamma$.\\

A \emph{holomorphic modular form} of integer weight $k \geq 0$ with respect to a congruence subgroup $\Gamma$ is a holomorphic function on the upper half-plane $\mathbb{H}$ which satisfies 
\begin{enumerate}
    \item[(i)] $f \big|_k \gamma = f$ for all $\gamma \in \Gamma$
    \item[(ii)] $f \big|_k \alpha$ is holomorphic at $i\infty$ for all $\alpha \in SL_2(\mathbb{Z})$.
\end{enumerate}
The function $f$ is said to be a \emph{weakly holomorphic modular form} if $f \big|_k \alpha$ is allowed to have poles at $i\infty$, that is, $f$ is \emph{meromorphic} at the cusps. More generally, $f$ is called a \emph{meromorphic modular form} if it is meromorphic on $\mathbb{H}$ and also at cusps. \\

We will say that \emph{$f$ has algebraic Fourier coefficients} if the Fourier coefficients of $f \big|_k \alpha$ at $i \infty$ are algebraic numbers. Denote the space of all holomorphic, weakly holomorphic and meromorphic modular forms, with algebraic Fourier coefficients by $M_{k,\overline{\mathbb{Q}}}(\Gamma)$, $M_{k,\overline{\mathbb{Q}}}^w(\Gamma)$ and $M_{k,\overline{\mathbb{Q}}}^m(\Gamma)$ respectively. Clearly, $M_{k,\overline{\mathbb{Q}}}(\Gamma)\subset M_{k,\overline{\mathbb{Q}}}^w(\Gamma)\subset M_{k,\overline{\mathbb{Q}}}^m(\Gamma)$. \\

For even integer $k\geq 2$, define the normalized Eisenstein series of weight $k$ for $S L_2(\mathbb{Z})$ by
$$
E_{k}(\tau) = 1-\frac{2 k}{B_{k}} \sum_{n=1}^{\infty} \sigma_{k-1}(n) q^{n}, \text{   where }q=e^{2\pi i\tau} \text{ and } \sigma_s(n)=\sum_{\substack{d|n\\ d>0}}d^s.
$$ 
Here, $B_k$ is the $k$-th Bernoulli number. Define
\begin{equation*}
    \Delta(\tau):=\frac{E_4(\tau)^3-E_6(\tau)^2}{1728}.
\end{equation*}
For $k\geq 4$, the function $E_k\in M_{k,\overline{\mathbb{Q}}}(S L_2(\mathbb{Z}))$ and $\Delta \in M_{12,\overline{\mathbb{Q}}}(S L_2(\mathbb{Z}))$. But $E_2(\tau)$ is not a modular form (see \cite{murty-dewar-graves}, Chapter $5$), as
\begin{equation*}
    E_2 \left( \frac{-1}{\tau} \right) = {\tau}^2 E_2(\tau) + \frac{6}{i\pi} \tau.
\end{equation*}
This motivates the definition of a quasi-modular form, of which there are several equivalent formulations. We use the following characterization, which was established for holomorphic quasi-modular forms in \cite[Proposition 20]{1-2-3} and meromorphic quasi-modular forms in \cite[Theorem 4.2]{wang-zhang}.
\begin{theorem}\label{quasistructure}
Every meromorphic quasi-modular form for a congruence subgroup $\Gamma$ is a polynomial in $E_2$ with modular coefficients. More precisely, if $\widetilde{f}$ is a meromorphic quasi-modular form of weight $k$ and depth $p$ with respect to $\Gamma$, then $\widetilde{f}$ can be uniquely written as $\widetilde{f} = \sum_{r=0}^p  f_r\, E_2^r$, where $f_r$ is a meromorphic modular form with respect to $\Gamma$ of weight $k-2r$ for all $0\leq r\leq p$ and $f_p\neq 0$. 
\end{theorem}
A quasi-modular form $\widetilde{f}$ is said to have algebraic Fourier coefficients if all the modular coefficients in the above expression of $f$ have algebraic Fourier coefficients. \\

\subsection{The Weierstrass $\wp$-function}
 Let $L=\omega_1 \mathbb{Z} \,  \oplus \, \omega_2 \mathbb{Z}$ be a two-dimensional lattice in $\mathbb{C}$ with $\omega_1/\omega_2 \in \mathbb{H}$. The Weierstrass $\wp$-function associated to $L$, given by
\begin{equation}\label{wpfunction}
\wp (z)=\frac{1}{z^2}+\sum_{\substack{\omega\in L\\
\omega\neq 0}}\left(\frac{1}{(z-\omega)^2}-\frac{1}{\omega^2} \right), \quad \text{for }z\in \mathbb{C}\setminus L,
\end{equation}
defines a meromorphic function on $\mathbb{C}$. It satisfies the differential equation
\begin{equation} \label{dell}
 \wp'(z)^2=4\wp(z)^3-g_2(L)\wp(z)-g_3(L),   
\end{equation}
where 
\begin{equation*}
    g_2(L)=60\, \sum_{\substack{\omega \in L\\ \omega\neq 0}} \frac{1}{\omega^{4}} \qquad \text{and} \qquad g_3(L)=140\, \sum_{\substack{\omega \in L\\ \omega\neq 0}} \frac{1}{\omega^{6}}.
\end{equation*}
For the lattice $L_\tau=\tau\,\mathbb{Z}\oplus \mathbb{Z}$ with $\tau\in\mathbb{H}$, 
\begin{equation*}
    g_2(L_\tau)=\frac{4\pi^4}{3}E_4(\tau) \qquad \text{and} \qquad g_3(L_\tau)=\frac{8\pi^6}{27}E_6(\tau).
\end{equation*}
Define the discriminant of a lattice, $$\Delta_0(L):=g_2(L)^3-27g_3(L)^2\neq 0, \text{ for any two dimensional lattice }L.$$ In particular, we have $\Delta_0(L_\tau)=(2\pi)^{12}\,\Delta(\tau)$ for all $\tau\in\mathbb{H}$.\\

The Weierstrass zeta-function associated to $L=\omega_1\,\mathbb{Z}\oplus \omega_2\,\mathbb{Z}$ is defined by
$$
\zeta_{L}(z)=\frac{1}{z}+\sum_{\substack{\omega\in L\\ \omega\neq 0}}\left(\frac{1}{z-\omega}+\frac{1}{\omega}+\frac{z}{\omega^2} \right)\quad \text{for } z\in \mathbb{C}\setminus L.
$$
Note that $\zeta_{L}'(z)=-\wp_L(z)$ is a periodic function with each point of $L$ as a period. Hence, the functions $$\zeta_{L}(z+\omega_1)-\zeta_{L}(z) \text{  and  } \zeta_{L}(z+\omega_2)-\zeta_{L}(z)$$ are constants. These constants are denoted by $\eta_1(L)$ and $ \eta_2(L)$ respectively and are called \emph{quasi-periods}. Moreover, for $\omega_1 = \tau$ and $\omega_2 = 1$, it is known that 
\begin{equation}\label{quasi-period}
\eta_2(L_\tau)=G_2(\tau)=\frac{\pi^2}{3}\, E_2(\tau) \text{ 
  for all }\tau\in\mathbb{H}. 
\end{equation}
The reader is referred to (\cite{lang}, chapter $18.3$) for proof of these results.\\

The Weierstrass $\wp$-function satisfies the following addition formula (\cite{r-m}, chapter $11$).
 \begin{theorem} \label{addf}
 For $z_1, z_2\in \mathbb{C}$ such that $z_1\pm z_2\notin L$ we have
$$
\wp\left(z_1+z_2\right)=-\wp\left(z_1\right)-\wp\left(z_2\right)+\frac{1}{4}\left(\frac{\wp^{\prime}\left(z_1\right)-\wp^{\prime}\left(z_2\right)}{\wp\left(z_1\right)-\wp\left(z_2\right)}\right)^2 .
$$
\end{theorem}
Using the addition formula of the Weierstrass $\wp$-function, Schneider \cite{schneider.} proved that 
\begin{theorem}\label{tranw}
    Let $L=\omega_1\,\mathbb{Z}\oplus\omega_2\,\mathbb{Z}$ be a lattice such that both $g_2(L)$, $g_3(L)$ are algebraic. If $\alpha$ is an algebraic number with $\alpha\notin L$, then $\wp(\alpha)$ is transcendental.
\end{theorem}
The addition formula also implies the following important proposition (see \cite{r-m}). Two generators $\omega_1$ and $\omega_2$ of a lattice $L$ are said to be primitive if both have minimal absolute value among all generators of $L$.
\begin{prop}\label{tranp}
Let $L=\omega_1\,\mathbb{Z}\oplus\omega_2\,\mathbb{Z}$ be such that both $g_2(L)$ and $g_3(L)$ are algebraic. Assume that $\omega_1$ and $\omega_2$ are primitive generators of $L$. Then, for any natural number $n>1$, the numbers $\wp(\frac{\omega_1}{n})$ and $\wp(\frac{\omega_2}{n})$ are algebraic. Moreover, any non-zero period of $L$ is necessarily transcendental.\\
\end{prop}

\subsection{Modular functions}

A meromorphic function $g$ on $\mathbb{H}$ is said to be a modular function if it satisfies
\begin{equation}
g\left(\frac{a\,\tau+b}{c\,\tau+d}\right)=g(\tau)\quad\text{for all }
\gamma=\left(\begin{smallmatrix}
   a  & b \\
   c  & d
\end{smallmatrix} \right)\in \Gamma,
\end{equation}
that is, $g \big|_0 \gamma = g$, and is also meromorphic at all the cusps. In particular, we call a modular function on the congruence subgroup $\Gamma(N)$ to a be \emph{modular function of level $N$}. Note that if $g$ is a modular function with respect to $\Gamma$, which is of level $N$, then $g$ is a modular function of level $N$.\\

An example of modular function of level one is given by 
$$
j(\tau):=\frac{E_4^3(\tau)}{\Delta(\tau)} \text{ for all } \tau\in\mathbb{H},
$$
which has the following Fourier expansion at $i\infty$
$$
j(\tau)=\frac{1}{q}+744+196884q+\cdots,$$
where $q=e^{2\pi i\tau}$. \\

More specifically, for any modular function (or form) $f$, define 
\begin{equation*}
    \mathbb{Q}(f) := \mathbb{Q} \left( \left\{ \text{Fourier coefficients of $f$ at all $\Gamma$-in-equivalent cusps} \right\} \right)
\end{equation*}
The $j$-function is a canonical example of a modular function of level one, and governs properties of modular functions of higher levels as well. This is made precise in the following series of propositions.
\begin{prop} \label{level1}
Let $g$ be a non-constant modular function of level one. Let $\mathcal{F}$ denote the standard fundamental domain for the action of $S L_2(\mathbb{Z})$ on $\mathbb{H}$. Suppose that the poles of $g$ in $\mathcal{F}$ are $\tau_1,\tau_2,\cdots,\tau_m$. Then $g(\tau)$ is a rational function in $j(\tau)$ with coefficients in the field $\mathbb{Q}(g) \left( j(\tau_1),\cdots,j(\tau_m) \right)$.
\end{prop}
\begin{proof}
By the compactness of $\mathcal{F} \cup i\infty$, we know that $g$ has only finitely many poles in $\mathcal{F}$. Consider the function 
$$
h(\tau):=\left(\prod_{j=1}^m \left(j(\tau)-j(\tau_j)\right)^{-\operatorname{ord}_{\tau_j}(g)}\right)g(\tau),
$$
where $\operatorname{ord}_{\tau_j}(g) = - $ order of the pole of $g$ at $\tau_j$.  Then $h$ is holomorphic on $\mathbb{H}$.\\

Suppose $h$ has pole of order $M \geq 1$ at $i\infty$. Then the Fourier expansion of
$h(\tau)$ at $i\infty$ has the form
$$
h(\tau)=\sum_{n=-M}^\infty c_n q^n, \, \text{ where } c_{-M}\neq 0.
$$
Note that the modular function $h(\tau)-c_{-M}\,j(\tau)^M$ is holomorphic on $\mathbb{H}$ and its Fourier expansion starts with at most a polar term of order $M-1$. Iterating this process, we can subtract a polynomial in $j(\tau)$ to get a holomorphic modular function that vanishes at $i\infty$, and hence is identically zero. Thus, $h(\tau)$ is a polynomial in $j(\tau)$ over $\mathbb{Q}(g)$, and so $g(\tau)$ is a rational function of $j(\tau)$ over $\mathbb{Q}(g)(j(\tau_1),\cdots,j(\tau_m))$.
\end{proof}

One can also conclude the following important fact from the above proof.
\begin{corollary}\label{rationalfunctions}
If $g$ is a modular function of level one which is holomorphic on $\mathbb{H}$ with a pole of order $M$ at $i\infty$, then $g(\tau)$ is a polynomial in $j(\tau)$ of degree $M$ with coefficients in $\mathbb{Q}(g)$.
\end{corollary}

The $j$-function is sufficient to `generate' all higher level modular functions as well. This is proved below.
\begin{theorem} \label{mdft}
Let $g$ be a modular function with respect to a congruence subgroup $\Gamma$ and let $\tau_1,\, \tau_2, \cdots, \, \tau_m$ be the poles of $g$ in $\mathcal{F}$. Set
$$
K:=\mathbb{Q}(g)\left(j(\tau_1), \, j(\tau_2), \cdots, \, j(\tau_m) \right).
$$
Then there exists a monic polynomial $P_g(X)\in K(j)[X]$ such that $P_g(g)=0$. 
\end{theorem}
\begin{proof}
Let $[S L_2(\mathbb{Z}):\Gamma]=r$ and $\{ \gamma_1=I,\gamma_2,\cdots,\gamma_r\}$ be a complete set of right coset representatives so that $$
S L_2(\mathbb{Z})=\bigsqcup_{i=1}^r\Gamma \gamma_i.
$$
For all $1\leq i \leq r$, define the functions $$g_i(\tau):=g(\gamma_i\tau) \quad \text{for } \tau \in \mathbb{H}.
$$ Each $g_i$ is independent of the choice of coset representatives as $g$ is $\Gamma$-invariant. Moreover, each $g_i$ is a modular function with respect to $\Gamma$ and the Fourier expansion of $g_i$ at $i \infty$ is precisely the expansion of $g$ at the cusp $\gamma_i(i\infty)$. For any $\gamma\in S L_2(\mathbb{Z})$, we have $g_i(\gamma \tau)=g(\gamma_i\gamma\tau)=g_j(\tau)$ for some $j$ with $1 \leq j \leq r$ such that $\gamma_i\gamma\in \Gamma\gamma_j$, i.e., the set $\{g_1,g_2,...,g_r\}$ gets permuted under the action of $S L_2(\mathbb{Z}).$ This observation implies that any elementary symmetric polynomial in $g_1,g_2,\cdots,g_r$ is a modular function of level one, and is in $K(j)$ by Theorem \ref{level1}. Note that the polynomial $$
P(X)=\prod_{i=1}^r(X-g_i)$$is satisfied by $g$ as $g=g_1$ and has coefficients that are elementary symmetric polynomials in $g_1,g_2,\cdots,g_r$. This proves the theorem.
\end{proof}

\begin{corollary}\label{useful}
Let $g$ be a modular function with respect to a congruence subgroup $\Gamma$ which is holomorphic on $\mathbb{H}$ and $K = \mathbb{Q}(g)$. Then the monic polynomial $P_g(X)$ satisfied by $g$ has coefficients in $K[j][X]$.
\end{corollary}
\begin{proof}
The coefficients of $P_g(X)$ constructed in the proof above are modular functions of level one, holomorphic on $\mathbb{H}$ with Fourier coefficients in $K$. The result now follows from Corollary \ref{rationalfunctions}. 
\end{proof}

\begin{remark*}
An important point to note here is that if $g$ is a modular function with respect $\Gamma$ that has algebraic Fourier coefficients, Theorem \ref{mdft} does not guarantee that the coefficients of $P_g$ are algebraic, unless $g$ is holomorphic in $\mathbb{H}$. 
\end{remark*}

To that effect, we now study the structure of the field of modular functions of a fixed level $N>1$. For any field $\mathbb{F} \subseteq \mathbb{C}$, let
\begin{equation*}
    \mathcal{F}_{N, \mathbb{F}} := \left\{ \text{Modular functions of level $N$ whose Fourier coefficients at $i \infty$ are in $\mathbb{F}$} \right\}.
\end{equation*}
Theorem \ref{mdft} implies that $\mathcal{F}_{N,\mathbb{C}}$ is an algebraic extension of $\mathbb{C}(j)$.

Following \cite[Chapter 6]{shimura}, we consider explicit modular function of level $N$ whose Fourier coefficients have good rationality properties. Let $\boldsymbol{a}=(a_1,a_2)\in \frac{1}{N}\mathbb{Z}^2 \,\setminus \,\mathbb{Z}^2$, consider the function
$$
f_{\boldsymbol{a}}(\tau):=\frac{g_2(L_\tau)\,g_3(L_\tau)}{\Delta_0(L_\tau)}\, \wp_{L_\tau} (a_1\tau+a_2),
$$
which is holomorphic on $\mathbb{H}$. The following properties can be checked routinely, and we leave the proof to the reader.
\begin{prop}
Let $S$ denotes the set $\{\left(\frac{r}{N},\frac{s}{N}\right): 0\leq r,s\leq N-1 \text{ and }(r,s)\neq(0,0)\}$. Let $\boldsymbol{a}$, $\boldsymbol{b} \in \frac{1}{N}\mathbb{Z}^2 \,\setminus \,\mathbb{Z}^2$ and $f_{\boldsymbol{a}}$ be as defined above. For $\gamma = \big(\begin{smallmatrix}
p & q\\
r & s
\end{smallmatrix}\big)$ and $\boldsymbol{a} = (a_1,a_2)$, let ${\boldsymbol{a\gamma}=(pa_1+ra_2, \, qa_1+sa_2)}$. Then
\begin{enumerate}
    \item[(a)] $f_{\boldsymbol{a}} = f_{\boldsymbol{b}}$ if and only if $\boldsymbol{a} \equiv \boldsymbol{b} \bmod \mathbb{Z}^2 $,
    \item[(b)] for $\gamma \in SL_2(\mathbb{Z})$, $f_{\boldsymbol{a}} (\gamma \tau) = f_{\boldsymbol{a\gamma}}(\tau)$ for all $\tau \in \mathbb{H}$.
\end{enumerate}
Therefore, all elements in $\{f_{\boldsymbol{a}}: \boldsymbol{a} \in S\}$ satisfy modularity with respect to $\Gamma(N)$.
\end{prop}

In order to conclude that the functions $f_{\boldsymbol{a}}$ are modular functions, we study their Fourier expansions at the cusp. But first, we need to understand the behaviour of the Weierstrass $\wp$-function.
\begin{lemma*}
Fix $\tau \in \mathbb{H}$ and denote $\wp_{\tau} = \wp_{L_{\tau}}$. For $z \not\in L_{\tau}$, we have
\begin{equation*}\label{weifou}
    (2\pi i)^{-2}\wp_\tau(z)=\frac{1}{12}+\frac{1}{\xi^{-1}-2+\xi}+\sum_{n=1}^\infty c_n q^n, 
\end{equation*}
where 
$$
\begin{aligned}
  \xi=e^{2\pi i z},\, q=e^{2\pi i\tau} \text{ and }
    c_n=\sum_{d|n}d\left(\xi^{-d}-2+\xi^d\right)\ \forall\ n\geq 1.
\end{aligned}
$$
\end{lemma*}
\begin{proof}
    From the definition of the $\wp$-function, we have 
$$
\begin{aligned}
\wp_\tau(z)&=\frac{1}{z^2}+\sum_{\substack{(m,n)\in \mathbb{Z}^2\\ (m,n)\neq (0,0)}}\left[\frac{1}{(z-m\tau-n)^2}-\frac{1}{(m\tau+n)^2} \right]\\
&=\frac{1}{z^2}+\sum_{n\in \mathbb{Z}\setminus \{0\}}\frac{1}{(z+n)^2}-\sum_{n\in \mathbb{Z}\setminus \{0\}}\frac{1}{n^2}+\sum_{m\neq 0}\sum_{n\in \mathbb{Z}}\frac{1}{(z-m\tau-n)^2}-\sum_{m\neq 0}\sum_{n\in \mathbb{Z}}\frac{1}{(m\tau+n)^2}\\
&=-2\zeta(2)+\sum_{n\in\mathbb{Z}}\frac{1}{(z+n)^2}-2\sum_{m=1}^\infty\sum_{n\in \mathbb{Z}}\frac{1}{(m\tau+n)^2}\\
& \quad \quad \quad \quad \quad \quad \quad \quad \quad \quad \quad +\sum_{m=1}^\infty\sum_{n\in \mathbb{Z}}\left[\frac{1}{(-z+m\tau+n)^2}+\frac{1}{(z+m\tau+n)^2} \right].
\end{aligned}
$$
Recall the Lipschitz summation formula, which states 
\begin{equation*} \label{lipc}
\sum_{n \in \mathbb{Z}} \frac{1}{(z+n)^{k}}=\frac{(-2 \pi i)^{k}}{(k-1) !} \sum_{n=1}^{\infty} n^{k-1} q^{n}.
\end{equation*}
For a proof, see \cite[Theorem 4.2.2]{murty-dewar-graves}. Applying this, we obtain
$$
\begin{aligned}
& \wp_\tau(z)=-\frac{\pi^2}{3}+(-2\pi i)^2\sum_{n=1}^\infty ne^{2\pi i n z}\\
& \quad \quad \quad \quad \quad\quad\quad\quad +(-2\pi i)^2\sum_{m=1}^\infty\sum_{n=1}^\infty n\left[e^{2\pi i n(-z+m\tau)}+e^{2\pi i n(z+m\tau)}-2e^{2\pi i m n \tau} \right] \\
&\Rightarrow (2\pi i)^{-2}\wp_\tau(z)=\frac{1}{12}+\sum_{n=1}^\infty ne^{2\pi i n z}+\sum_{m,n=1}^\infty nq^{mn}\left[e^{2\pi i n z}+e^{-2 \pi i n z}-2 \right]
  \end{aligned}
$$
This implies the result.
\end{proof}

Using the above expansion, the interpretation of $g_2(L_{\tau)}$, $g_3(L_{\tau})$ and $\Delta_{0}(L_{\tau})$ in terms of Eisenstein series, we get for $\boldsymbol{a}= \left( \frac{r}{N}, \frac{s}{N} \right) \in S$,
\begin{align}
& f_{\boldsymbol{a}}(\tau) \nonumber \\
& =\frac{-1}{2592}\cdot\frac{E_4(\tau)E_6(\tau)}{\Delta(\tau)}\left[\frac{1}{12} -\sum_{n=1}^\infty ne^{\frac{2\pi i n s}{N}} q_N^{n r}+\sum_{m,n=1}^\infty n q^{m n}\left[e^{\frac{2\pi i n s}{N}}q_N^{n r}+e^{\frac{-2\pi i n s}{N}}q_N^{-n r}-2 \right]\right]. \label{foufa}
\end{align}
Here $q_N = e^{2 \pi i \tau/N}$. Since the Fourier series of $\Delta(\tau)$ begins with $q=q_N^N$, the Fourier expansion of $f_{\boldsymbol{a}}(\tau)$  begins with a rational multiple of $q_N^{-N}$. Thus, $f_{\boldsymbol{a}}$ has pole of order $N$ at $i\infty$ for all $\boldsymbol{a}\in S$. If $s$ is any other cusp, then there exists $\gamma\in S L_2(\mathbb{Z})$ such that $s=\gamma(i\infty)$. Since $f_{\boldsymbol{a}}(\gamma\tau)=f_{\boldsymbol{a\gamma}}(\tau)$, which has pole of order $N$ at $i\infty$, we conclude that $f_{\boldsymbol{a}}$ is meromorphic at all cusps, with pole of order $N$. Thus, $f_{\boldsymbol{a}}$ is a modular function of level $N$ for all $ \boldsymbol{a} \in S$. \\

This helps us to conclude the following.
\begin{theorem}\label{primiext}
For all $\boldsymbol{a} \in S$, the Fourier coefficients of $f_{\boldsymbol{a}}$ with respect to all cusps belong to $\mathbb{Q}(\mu_N)$, where $\mu_N=e^{2\pi i/N}$. 
\end{theorem}
\begin{proof}
   Recall that Fourier coefficients of $E_4(\tau), E_6(\tau)$ and $\Delta(\tau)$ are integers. Hence, from \eqref{foufa}, it follows that the Fourier coefficients of $f_{\boldsymbol{a}}$ at $i\infty$ lie in $\mathbb{Q}(\mu_N)$ for all $\boldsymbol{a}\in S$. If $s$ is any other cusp, then there exists $\gamma\in S L_2(\mathbb{Z})$ such $s=\gamma(i\infty)$. But $f_{\boldsymbol{a}}(\gamma\tau)=f_{\boldsymbol{a\gamma}}(\tau)$ also has Fourier coefficients in $\mathbb{Q}(\mu_N)$ with respect to $i\infty$. This completes the proof.
\end{proof}

These modular functions, together with the $j$-function serve to generate all modular functions of level $N$. That is,
\begin{theorem}\label{mdft..}
We have $\mathcal{F}_{N,\overline{\mathbb{Q}}}=\overline{\mathbb{Q}}\left( j, \, \{ f_{\boldsymbol{a}} |{\boldsymbol{a}}\in S\} \right)$ and $\mathcal{F}_{N,\overline{\mathbb{Q}}}$ is a finite Galois extension of $\overline{\mathbb{Q}}(j)$.
\end{theorem}
\begin{proof}

Let $E_{N,\overline{\mathbb{Q}}}:= \overline{\mathbb{Q}}\left(j, \, \{ f_{\boldsymbol{a}} | {\boldsymbol{a}}\in S\}\right)$. Then $E_{N, \overline{\mathbb{Q}}}$ is a Galois extension of $\overline{\mathbb{Q}}(j)$. Indeed, for each ${\boldsymbol{a}}\in S$, the modular function $f_{\boldsymbol{a}}$ is holomorphic on $\mathbb{H}$ with algebraic Fourier coefficients at all cusps. Hence, by Corollary \ref{useful} we get a polynomial $P(X)\in\overline{\mathbb{Q}}(j)[X]$, which is satisfied by $f_{\boldsymbol{a}}$. Thus, $E_{N,\overline{\mathbb{Q}}} \, / \, \overline{\mathbb{Q}}(j)$ is an algebraic extension. Moreover, this extension is normal because each conjugate of $f_{\boldsymbol{a}}$ is of the form $f_{\boldsymbol{a\gamma_i}}$, which is equal to $f_{\boldsymbol{b}}$ for some $\boldsymbol{b}\in S$. Since $E_{N,\overline{\mathbb{Q}}}$ is obtained from $\overline{\mathbb{Q}}(j)$ by adjoining finite number of elements, $E_{N,\overline{\mathbb{Q}}} \, / \, \overline{\mathbb{Q}}(j)$ is a finite Galois extension.\\

Now we show that $E_{N,\overline{\mathbb{Q}}}= \mathcal{F}_{N,\overline{\mathbb{Q}}}$. To begin with, observe that $\mathbb{C}$ and $\mathcal{F}_{N,\overline{\mathbb{Q}}}$ are linearly disjoint over $\overline{\mathbb{Q}}$. This can be seen as follows. Suppose $\{c_1,c_2,\ldots,c_r\}\subset \mathbb{C}$ is an arbitrary $\overline{\mathbb{Q}}$-linearly independent subset of $\mathbb{C}$. If there exists $g_i \in \mathcal{F}_{N,\overline{\mathbb{Q}}}$ such that $\sum_{i=1}^r g_i(\tau) \, c_i = 0$ for all $\tau \in \mathbb{H}$ with
$$g_i(\tau)=\sum_n d_{i n} \, \, q_N^n, \qquad d_{in} \in \overline{\mathbb{Q}} \quad \text{for }0\leq i\leq r,$$
then $$\sum_{i=1}^rc_i\,d_{i n}=0\quad \text{ for all } n.$$ The linear independence of $c_i$ over $\overline{\mathbb{Q}}$ implies that all $d_{i n}=0$ for all $0\leq i\leq r$, and hence, $f_1=f_2=\cdots=f_r=0$. Thus, we have $E_{N,\overline{\mathbb{Q}}}\subseteq \mathcal{F}_{N,\overline{\mathbb{Q}}} \subseteq \mathbb{C}E_{N,\overline{\mathbb{Q}}}.$ Suppose there exists $f\in \mathcal{F}_{N,\overline{\mathbb{Q}}}\setminus E_{N,\overline{\mathbb{Q}}}$. Since $f \in \mathbb{C}E_{N,\overline{\mathbb{Q}}}$, we get a $\overline{\mathbb{Q}}$-linearly independent subset $\{f_1,f_2,\ldots,f_m\}\subseteq E_{N,\overline{\mathbb{Q}}}$ such that 
\begin{equation}\label{rep.}
  f=\sum_{i=1}^m \alpha_i f_i, \quad \alpha_i\in \mathbb{C},  
\end{equation}
with at least one of the $\alpha_i \in \mathbb{C} \setminus \overline{\mathbb{Q}}$. Since $E_{N,\overline{\mathbb{Q}}}\subseteq \mathcal{F}_{N,\overline{\mathbb{Q}}}$, the set $\{f, \,f_1,f_2,\ldots,f_m\}$ is $\overline{\mathbb{Q}}$-linearly independent subset of $\mathcal{F}_{N,\overline{\mathbb{Q}}}$ and hence, $\mathbb{C}$-linearly independent. This contradicts \eqref{rep.}. Therefore, $E_{N,\overline{\mathbb{Q}}}= \mathcal{F}_{N,\overline{\mathbb{Q}}}$ and $\mathcal{F}_{N,\overline{\mathbb{Q}}} \, / \, \overline{\mathbb{Q}}(j)$ is a finite Galois extension.
\end{proof}

Thus, Theorem \ref{primiext} and Theorem \ref{mdft..} imply the following crucial fact.
\begin{corollary}\label{cor-main-mdf}
Let $g \in \mathcal{F}_{N,\overline{\mathbb{Q}}}$. Then the Fourier coefficients of $g$ with respect to all cusps are algebraic numbers and $g$ satisfies a polynomial over $\overline{\mathbb{Q}}(j)$. 
\end{corollary}

\subsection{Algebraic independence of modular, quasi-modular forms and modular functions} 

The aim of the discussion below is to establish the algebraic independence (as functions) of the three functions arising from modular considerations - modular functions, modular forms and quasi-modular forms. For the basic theory of quasi-modular forms, we refer the reader to \cite{royer}. To begin with, we prove the following lemma.
\begin{lemma}\label{indep}
A sum of meromorphic quasi-modular forms of distinct weights is not identically zero, unless, each of these quasi-modular forms is identically zero.
\end{lemma}
\begin{proof}
Suppose $f_1, f_2,\ldots,f_r$ are non-identically zero meromorphic quasi-modular forms of weights $k_1<k_2<\cdots<k_r$ with respect to the congruence subgroups $\Gamma_1,\Gamma_2,\ldots,\Gamma_r$ of level $N_1,N_2,\ldots,N_r$ respectively. Let $p$ be the greatest depth of these quasi-modular forms. Suppose that for each $i\in\{1,2,\ldots,r\}$, they have the following transformation formulae:
$$
\begin{aligned}
&f_i\big|_{k_i}\left(\begin{smallmatrix}
a & b\\
c & d
\end{smallmatrix}\right)(\tau)=\sum_{j=0}^p f_{i,j}(\tau)\left(\frac{c}{c\tau+d}\right)^j\\
\text{i.e.,}\quad & f_i\left(\frac{a\tau+b}{c\tau+d}\right)=\sum_{j=0}^p f_{i,j}(\tau)\,c^j\, (c\tau+d)^{k_i-j}
 \end{aligned}
$$
for every $\left(\begin{smallmatrix}
a & b\\
c & d
\end{smallmatrix}\right)\in \Gamma_i$ and $\tau\in\mathbb{H}$. Here the functions $f_{i,0},f_{i,1},\ldots, f_{i,p}$ are the components of $f_i$, and in particular, $f_{i,0}=f_i$ for all $1\leq i\leq r$. Then we consider $N:=\prod_{i=1}^r N_i$ and $S:=\{bN^2+1: b\in\mathbb{N}\}$ such that for each $b'=Nb^2+1\in S$, the matrix $\left(\begin{smallmatrix}
1 & bN\\
N & b'
\end{smallmatrix}\right)\in \cap_{i=1}^r\Gamma_i$. For all such matrices we have 
$$
f_i\left(\frac{\tau+bN}{N\tau+b'}\right)=\sum_{j=0}^p f_{i,j}(\tau)\,N^j\, \left(N\tau+b'\right)^{k_i-j} \quad \text{  for all } 1\leq i\leq r.
$$
Suppose that $\sum_{i=1}^r f_i=0$. Then for $\tau\in\mathbb{H}$ and $b'\in S$, we have 
$$
\sum_{i=1}^r f_i\left(\frac{\tau+bN}{N\tau+b'}\right)=0.
$$
 From the above transformation formula for each $f_i$, we obtain
\begin{equation}\label{1}
\sum_{i=1}^r\sum_{j=0}^p f_{i,j}(\tau)\, N^j\, (N\tau+b')^{k_i-j}=0.
 \end{equation}
For a fixed $\tau\in\mathbb{H}$, multiplying \eqref{1} by $(N\tau+b')^{p'}$, where $p'=2\cdot\max\{p,|k_1|,\ldots,|k_r|\} $,  we get $$
\begin{aligned}
0&=\sum_{i=1}^r\sum_{j=0}^p f_{i,j}(\tau)\, N^j\, (N\tau+b')^{k_i+p'-j}\\
&=\sum_{i=1}^r N^{k_i+p'}\sum_{j=0}^pf_{i,j}(\tau)\left(\tau+\frac{b'}{N}\right)^{k_i+p'-j},
\end{aligned}
$$
which holds for each $b'\in S$. This shows that for any fixed $\tau\in \mathbb{H}$, the polynomial 
$$
P(X)=\sum_{i=1}^r N^{k_i+p'}\sum_{j=0}^pf_{i,j}(\tau)\left(\tau+X\right)^{k_i+p'-j}\in \mathbb{C}[X]
$$
has infinitely many roots, and hence $P(X)$ is identically zero. Note that even if any of the $k_j$'s are negative, the maximum exponent of $X$ occurring in $P(X)$ is $p'+ k_r$. Thus, the leading coefficient of $P(X)$ is $N^{k_r+p'}\, f_{r,0}(\tau)$. Since $P(X)=0$, we get $f_{r,0}(\tau)=f_r(\tau)=0$. This is true for any $\tau\in\mathbb{H}$. Hence, $f_r=0$, which is a contradiction, proving the lemma.  
\end{proof}

We are now ready to prove the main theorem in this context.
\begin{theorem}\label{algebraic-indep-fn-thm}
Let $f$ a non-constant meromorphic modular form of non-zero weight. Let $g$ and $\widetilde{f}$ be a modular function and a meromorphic quasi-modular form of positive depth and non-zero weight, respectively. Assume that $f$, $g$ and $\widetilde{f}$ are of arbitrary level. Then the functions $f$, $g$ and $\widetilde{f}$ are algebraically independent. 
\end{theorem}
\begin{proof}
Let $f$ be of weight $m$ and level $N$, $g$ be of level $M$ and $\widetilde{f}$ be of weight $n$ and level $\widetilde{N}$. Suppose that $P\in \mathbb{C}[X,Y,Z]$ is such that $P(f,g,\widetilde{f})=0$. Recall that a product of a modular form of weight $m$ and level $N$ and a quasimodular form of weight $n$ and level $\widetilde{N}$ is a quasimodular form of weight $mn$ and level $N\widetilde{N}$. Denote $\mathbb{N}_0 = \mathbb{Z}_{\geq 0}$. By grouping the terms of $P(f,g,\widetilde{f})$ of the same weight, we can rewrite it as
\begin{equation}\label{2}
P(f,g,\widetilde{f})=\sum_{k=-K'}^K\,\sum_{\substack{{(r,s,t)\in\mathbb{N}_{0}^3}\\ {mr+nt=k}}} p_{r,s,t} \, \, f^r\,g^s\,{\widetilde{f}}^t , \qquad \text{with   } p_{r,s,t} \in \mathbb{C}.
 \end{equation} 
As $P(X,Y,Z)$ is a polynomial, $p_{r,s,t}=0$ for all but finitely many $r$, $s$ and $t$. As noted above, for each $k\in \{-K',\ldots,K\}$ the inner sum in \eqref{2} is a meromorphic quasi-modular form of weight $k$ and level $NM\widetilde{N}$. Thus, Lemma \ref{indep} implies that
\begin{equation}\label{3}
\sum_{\substack{{(r,s,t)\in\mathbb{N}_0^3}\\ {mr+nt=k}}} p_{r,s,t} \,\, f^r\,g^s\,{\widetilde{f}}^t=0.
\end{equation}
 for each $-K'\leq k\leq K$. \\
 
\noindent Fix a $k\in\{-K',\ldots,K\} $ and consider the corresponding inner sum from \eqref{3}. If there exists an integer $t_0\neq 0$ such that $p_{r,s,t_0}\neq 0$ for at least one tuple $(r,s,t_0)\in\mathbb{N}_0^3$ satisfying $mr+nt_0=k$, then the term on the left in \eqref{3}
is a meromorphic quasi-modular form of depth at least $t_0$. But the uniqueness of depth implies that the term on the left in \eqref{3} is of depth $0$. Hence, the coefficient $p_{r,s,t}=0$ when $t\neq 0$. If $t=0$, then $mr=k$, or in other words, $r = k/m \in \mathbb{Z}$. Thus, for all $0 \neq k\in\{-K',\ldots, K\}$ such that $\frac{k}{m}\in\mathbb{Z}$, the relation in \eqref{3} takes the form  
\begin{equation}\label{5}
\sum_{s \in \mathbb{N}_0} p_{\frac{k}{m},s,0} \,\,g^s = 0,
\end{equation}
as $f$ is not identically zero. If $k=0$, then clearly, $r=0$ and we get
\begin{equation}\label{4}
    \sum_{s \in \mathbb{N}_0} p_{0,s,0} \, \, g^s=0.
\end{equation}
Since $g$ is a non-constant modular function, $g$ must have either a zero or a pole in $\mathbb{H}$. Say that the order of $g$ at $\tau_0$ is $b \neq 0$. Then any function $\sum_{r=0}^R c_r\, g^r$ with not all $c_r = 0$ will have order $bR \neq 0$ at $\tau_0$. Thus, any combination of the form  $\sum_{r=0}^R c_r\, g^r$ cannot be identically zero. Thus, relations \eqref{5} and \eqref{4} imply that $p_{0,s,0}=p_{\frac{k}{m},s,0}=0$ for all $k$ and $s$. Hence $P=0$ and the theorem is proved.
\end{proof}

\bigskip

\section{\bf Proof of Main Results}
\bigskip

\subsection{Values of modular functions when $j(\tau) \not\in \overline{\mathbb{Q}}$}
Before proceeding with the proof of the main theorems, we establish the following intermediary observation.
\begin{lemma} \label{mdfz}
    Let $g\in  \mathcal{F}_{N,\overline{\mathbb{Q}}}$. If $\tau\in \mathbb{H}$ such that $j(\tau)\notin\overline{\mathbb{Q}}$, then $g(\tau)$ is non-zero.
\end{lemma}
\begin{proof}
By Theorem \ref{mdft..}, we know that $g$ satisfies a non-trivial polynomial over $\overline{\mathbb{Q}}(j)$. On clearing denominators, we can assume that $g$ satisfies the irreducible polynomial
$$
P(X)= \sum_{r=0}^m c_r(j) \, X^r =\sum_{r=0}^m \bigg( \sum_{s=0}^{d_r}c_{r s}\,j^s \bigg) \,X^r \in \overline{\mathbb{Q}}[j][X], \quad m=\deg P(X), \, \, d_r=\deg c_r(j),
$$
Note that $c_0(j) \in \overline{\mathbb{Q}}[j]$ is not the identically zero polynomial. Now suppose $\tau \in \mathbb{H}$ is such that $j(\tau)$ is transcendental. If $g(\tau)$ is zero, then we get that $c_0(j(\tau))=0$ which contradicts the transcendence of $j(\tau)$. This completes the proof. 
\end{proof}

\begin{proof}[\bf Proof of Theorem \ref{mdf-values-thm}(a)]
If $g(\tau)=0$, then it follows from Lemma \ref{mdfz} that $j(\tau)\in \overline{\mathbb{Q}}$. If $\tau$ is a pole of $g$, we consider the modular function
$$
h(z):=\frac{1}{g(z)}\in \mathcal{F}_{N,\overline{\mathbb{Q}}},
$$
which vanishes at $z=\tau$, and so $j(\tau)\in\overline{\mathbb{Q}}$ by Lemma \ref{mdfz}.    
\end{proof}

The proof of Theorem \ref{mdf-values-thm}(b) is established in two parts. We first show that if $j(\tau)$ is transcendental, then $g(\tau)$ is as well. The converse implication is proved later.
\begin{prop}\label{tranmf}
 Let $g\in \mathcal{F}_{N,\overline{\mathbb{Q}}}$ and $\tau\in \mathbb{H}$ be such that it is not a pole of $g$. If $j(\tau)$ is transcendental, then $g(\tau)$ is transcendental.
\end{prop}
\begin{proof}
Suppose that $g(\tau)\in\overline{\mathbb{Q}}$. Then the modular function $$h(z):=g(z)-g(\tau)$$ belongs to $ \mathcal{F}_{N,\overline{\mathbb{Q}}}$ and vanishes at $\tau$. This contradicts Lemma \ref{mdfz}. Therefore, $g(\tau)$ is transcendental.
\end{proof}

\subsection{Values of modular functions when $j(\tau)\in\overline{\mathbb{Q}}$}
From Theorem \ref{mdft..}, we have $$ \mathcal{F}_{N,\overline{\mathbb{Q}}}=  \overline{\mathbb{Q}}\left(\{j,f_{\boldsymbol{a}} : {\boldsymbol{a}}\in S\}\right), \text{ where }S=\left\{\left(\frac{r}{N},\frac{s}{N}\right): 0\leq r, s\leq N-1 \text{ and }(r,s)\neq (0,0)\right\}.$$
Hence, for any $g\in \mathcal{F}_{N,\overline{\mathbb{Q}}}$, the algebraic nature of $g(\tau)$ is determined by the numbers in the set $\{f_{\boldsymbol{a}}(\tau):{\boldsymbol{a}}\in S\}$.\\

For $N=1$, $\mathcal{F}_{1,\overline{\mathbb{Q}}} = \overline{\mathbb{Q}}(j)$. Hence, for any $g \in \mathcal{F}_{1,\overline{\mathbb{Q}}}$, $j(\tau) \in \overline{\mathbb{Q}}$ implies that $g(\tau) \in \overline{\mathbb{Q}}$. For $N>1$, we have to study the nature of the values $f_{\boldsymbol{a}}(\tau)$ for all ${\boldsymbol{a}}\in S$. We start by proving the following important lemma.
\begin{lemma}\label{omegatau}
 If $\tau\in \mathbb{H}$ such that $j(\tau)$ is algebraic, then there exists a unique transcendental number (up to an algebraic multiples) $\omega_\tau$ for which $g_2(\omega_\tau L_\tau)$ and $g_3(\omega_\tau L_\tau)$ are both algebraic numbers.   
\end{lemma}
\begin{proof}
Note that for any $\alpha \in\mathbb{C}^\times $,
\begin{equation}\label{omegatau.}
    j(\tau)=\frac{E_4(\tau)^3}{\Delta(\tau)}=\frac{1728\,  g_2(L_\tau)^3}{\Delta_0(L_\tau)}=\frac{1728\, g_2\left(\alpha L_\tau\right)^3}{\Delta_0\left(\alpha L_\tau\right)},
\end{equation}
because of the homogeneity properties of the $g_2$, $g_3$ functions. Here $L_\tau=\tau\,\mathbb{Z}\oplus\mathbb{Z}$ and $\alpha L_\tau= \alpha \, \tau\,\mathbb{Z}\oplus\alpha \, \mathbb{Z}$. If we choose $\omega_\tau$ such that $\omega_\tau^4=g_2(L_\tau)$, then $g_2(\omega_\tau L_\tau)=1$. Given that $j(\tau)$ is algebraic, from \eqref{omegatau.}, we get that $\Delta_0(\omega_\tau L_\tau)$ is also algebraic. Thus, the numbers $g_2(\omega_\tau L_\tau)$ and $g_3(\omega_\tau L_\tau)$ are both algebraic. Moreover, $\omega_\tau$ is a period of $\omega_\tau L_\tau$ and both $g_2(\omega_\tau L_\tau)$ and $g_3(\omega_\tau L_\tau)$ are algebraic. By Corollary \ref{tranp} we conclude that $\omega_\tau$ is transcendental.\\

To prove uniqueness, consider an arbitrary $\omega_\tau'\in\mathbb{C}^\times$ such that both $g_2(\omega_\tau' L_\tau)$ and $g_3(\omega_\tau' L_\tau)$ are algebraic. By homogeneity of $g_2$, we have 
$$
 \quad g_2(\omega_\tau' L_\tau)= (\omega_\tau')^{-4}\,g_2(L_\tau)=:\beta \in \overline{\mathbb{Q} }.
$$
Then, we get $g_2(L_\tau)=\beta \cdot (\omega_\tau')^{4}$ and hence $\omega_\tau^4=\beta \cdot (\omega_\tau')^{4}.$
This completes the proof that $\omega_\tau$ is unique up to algebraic multiples.
\end{proof}

Now recall that
$$
f_{\boldsymbol{a}}(\tau):=\frac{g_2(L_\tau)\,g_3(L_\tau)}{\Delta_0(L_\tau)}\,\wp_{L_\tau}(a_1\tau+a_2),
$$ 
where $L_\tau=\tau\,\mathbb{Z}\oplus\mathbb{Z}$ and ${\boldsymbol{a}}=(a_1,a_2)\in S$. Using homogeneity properties of the functions involved, we can rewrite this as follows:
$$
f_{\boldsymbol{a}}(\tau)=\frac{g_2(\omega_\tau L_\tau)\, g_3(\omega_\tau L_\tau)}{\Delta_0(\omega_\tau L_\tau)}\,\wp_{\omega_\tau L_\tau}(a_1\omega_\tau\tau+a_2\omega_\tau).
$$
\begin{lemma}\label{algfa}
    If $\tau\in\mathbb{H}$ is such that $j(\tau)\in\overline{\mathbb{Q}}$, then $f_{\boldsymbol{a}}(\tau)\in\overline{\mathbb{Q}}$ for all ${\boldsymbol{a}}\in S$.
\end{lemma}
\begin{proof}
    Let ${\boldsymbol{a}}=\left( \frac{r}{N},\frac{s}{N} \right)\in S$. Then we have
$$
f_{\boldsymbol{a}}(\tau)=\frac{g_2(\omega_\tau L_\tau)\, g_3(\omega_\tau L_\tau)}{\Delta_0(\omega_\tau L_\tau)}\,\wp_{\omega_\tau L_\tau}\left(\frac{r\,\omega_\tau\,\tau}{N}+\frac{s\,\omega_\tau}{N}\right),
$$
where $\omega_\tau$ is chosen as in Lemma \ref{omegatau} so that $g_2(\omega_\tau L_\tau),\, g_3(\omega_\tau L_\tau)\in\overline{\mathbb{Q}}$.\\

Recall that both $\wp_{\omega_\tau L_\tau}\left(\frac{\omega_\tau\, \tau}{N} \right)$ and $\wp_{\omega_\tau L_\tau}\left(\frac{\omega_\tau}{N} \right)$  are algebraic by Corollary \ref{tranp}. Moreover, by the addition formula for Weierstrass $\wp$-function in Theorem \ref{addf}, we get that $$\wp_{\omega_\tau L_\tau}\left(\frac{r\,\omega_\tau\, \tau}{N}+\frac{s\,\omega_\tau}{N}\right)\in\overline{\mathbb{Q}}$$ for all ${\boldsymbol{a}}=\left(\frac{r}{N},\frac{s}{N}\right)\in S$. Thus, the number $f_{\boldsymbol{a}}(\tau)$ is algebraic for all ${\boldsymbol{a}}\in S$.
\end{proof}
\begin{prop}\label{algemf}
Let $g\in \mathcal{F}_{N,\overline{\mathbb{{Q}}}}$ and $\tau\in\mathbb{H}$ such that $j(\tau)\in \overline{\mathbb{Q}}$. If $\tau$ is not a pole of $g$, then $g(\tau)$ is algebraic.
\end{prop}
\begin{proof}
 From Lemma \ref{algfa}, we have   $f_{\boldsymbol{a}}(\tau)\in\overline{\mathbb{Q}}$ for all ${\boldsymbol{a}}\in S$. Since 
 $$g \in \mathcal{F}_{N,\overline{\mathbb{Q}}}=  \overline{\mathbb{Q}}\left(\{j,f_{\boldsymbol{a}} |{\boldsymbol{a}}\in S\}\right)$$
 by Theorem \ref{mdft..}, $g(\tau)$ is also algebraic. 
\end{proof}
\bigskip

Propositions \ref{tranmf} and \ref{algemf} together complete the proof of Theorem \ref{mdf-values-thm}(b). We now prove the generalization of \cite[Theorem 2.4]{jeon-kang-kim}.
\begin{proof}[\bf Proof of Theorem \ref{zero-pole-value}]
Let $f$ be a meromorphic modular form of weight $k \in \mathbb{Z}$. Define 
   $$
h(\tau):=\frac{f(\tau)^{12}}{\Delta(\tau)^k}\quad \forall \tau\in\mathbb{H}.
   $$
The function $h \in \mathcal{F}_{N,\overline{\mathbb{Q}}}$ because the Fourier coefficients of $f$ and $\Delta$ are algebraic. Moreover, any $\tau\in\mathbb{H}$ is a zero or a pole of $h$ if and only if it is a zero or a pole of $f$ because $\Delta$ is non-vanishing and holomorphic on $\mathbb{H}$. Thus, if $\tau\in\mathbb{H}$ is a zero or pole of $f$, then $h$ is non-constant and by Theorem \ref{mdf-values-thm}(a), we deduce that $j(\tau)\in\overline{\mathbb{Q}}$, and by Theorem \ref{mdf-values-thm}(b), that either $g(z)$ has a pole at $z = \tau$ or that $g(\tau)\in\overline{\mathbb{Q}}$.
\end{proof}
\bigskip

\subsection{Values of modular forms} Recall that for any integer $k$ and $N$, we denote $M_{k,N,\overline{\mathbb{Q}}}$, $M_{k,N,\overline{\mathbb{Q}}}^w$ and $M_{k,N,\overline{\mathbb{Q}}}^m$ to be the set of holomorphic, weakly holomorphic and meromorphic modular forms respectively of level $N$ with algebraic Fourier coefficients at $i\infty$. By Corollary \ref{cor-main-mdf}, all elements of these sets have algebraic Fourier coefficients with respect to all cusps as well.\\

\begin{proof}[Proof of Theorem \ref{general-mero-values-thm}(a) ]
Let $f \in M^m_{k,N,\overline{\mathbb{Q}}}$. Consider the modular function 
\begin{equation*}
    g(\tau) := \frac{{f(\tau)}^{12}}{\Delta(\tau)^k}.
\end{equation*}
If $g$ is a non-zero constant, then $f(\tau) = c \Delta(\tau)^{k/12}$ for some $c \in \overline{\mathbb{Q}}^{\times}$ and all $\tau \in \mathbb{H}$. Since $1728 \Delta(\tau) = E_4^3(\tau) - E_6^2(\tau)$, and $E_4(\tau)$ is algebraically independent with $E_6(\tau)$ for $\tau \in \mathbb{H}$ such that $e^{2 \pi i \tau} \in \overline{\mathbb{Q}}$ by Nesterenko's theorem \ref{nest}, we deduce that $f(\tau)$ is transcendental. \\

Now suppose that $g$ is non-constant. Since $f \in M_{k,N,\overline{\mathbb{Q}}}$, $g \in \mathcal{F}_{N,\overline{\mathbb{Q}}}$. Note that for $\tau_0 \in \mathbb{H}$ such that $e^{2 \pi i \tau_0}$ is algebraic, $j(\tau_0)$ is transcendental by Theorem \ref{mahl}. Such a $\tau_0$ cannot be a pole of $f$. For if $\tau_0$ is a pole of $f$, then $\tau_0$ is a pole of $g$ and by Theorem \ref{mdf-values-thm}(a), $j(\tau_0)$ would be algebraic, leading to a contradiction. As $\tau_0$ is not a pole of $g$, Theorem \ref{mdf-values-thm}(b) implies that $g(\tau_0)$ is transcendental. Moreover, by Theorem \ref{mdft..}, $g(\tau_0)$ is algebraic over $\overline{\mathbb{Q}}(j(\tau_0))$. By Nesterenko's theorem, 
\begin{equation*}
    \Delta(\tau_0) = \frac{E_4^3(\tau_0) - E_6^2(\tau_0)}{1728} \qquad \text{ and } \qquad j(\tau_0) = \frac{E_4^3(\tau_0)}{\Delta(\tau_0)}
\end{equation*}
are algebraically independent. Therefore,
\begin{equation*}
    \text{trdeg}_{\,\mathbb{Q}} \, \, \overline{\mathbb{Q}} \left ( g(\tau_0), \, \Delta(\tau_0) \right) = \text{trdeg}_{\,\mathbb{Q}}\, \, \overline{\mathbb{Q}} \left ( j(\tau_0), \, \Delta(\tau_0) \right) = 2.
\end{equation*}
Hence, $f(\tau_0)^{12} = g(\tau_0) \, \Delta(\tau_0)^k$ is transcendental, proving the claim.
\end{proof}

We remark here that the proof of Theorem \ref{general-mero-values-thm}(a) only requires Nesterenko's theorem and the structure of modular functions of higher level. One also immediately deduces the following. 
\begin{prop}\label{algi-cor}
If $f\in M^m_{k,N,\overline{\mathbb{Q}}}$, then $f$ is algebraically dependent with $E_4$ and $E_6$ over $\overline{\mathbb{Q}}$. In particular, If $f\in M_{k,N,\overline{\mathbb{Q}}}^m$ and $\tau\in\mathbb{H}$ is not a pole of $f$, then $f(\tau)$ is algebraic over $\overline{\mathbb{Q}}(E_4(\tau),E_6(\tau))$. 
\end{prop}
\begin{proof}
Let $f \in M^m_{k,N,\overline{\mathbb{Q}}}$ and consider 
\begin{equation*}
    g (\tau) := \frac{f^{12}(\tau)}{\Delta^k(\tau)} \in \mathcal{F}_{N,\overline{\mathbb{Q}}}.
\end{equation*}
By Corollary \ref{cor-main-mdf}, there exists a polynomial $P(X) \in \overline{\mathbb{Q}}[j](X)$ such that $P(g) = 0$. More specifically,
$$\sum_{r=0}^m \sum_{s=0}^{d_r} c_{r,s} \, {j(\tau)}^s \,{g(\tau)}^r= 0 \qquad \text{for all} \qquad \tau \in \mathbb{H}.$$
Here $c_{r,s} \in \overline{\mathbb{Q}}$ for all $0 \leq s \leq d_r$ and $0 \leq r \leq m$. Multiplying by $\Delta^l(\tau)$ for any positive integer $l > km$ and substituting $j(\tau)$ and $\Delta(\tau)$ in terms of $E_4(\tau)$ and $E_6(\tau)$ gives
\begin{equation*} 
\sum_{r=0}^m \sum_{s=0}^{d_r} \sum_{t=0}^{l-k r-s} \frac{(-1)^t c_{r, s}}{1728^{l-k r-s}}\left(\begin{array}{c}
l-k\,r-s \\
t
\end{array}\right) f(\tau)^{12\, r} E_4(\tau)^{3(l-k\, r-t)} E_6(\tau)^{2\, t}=0.   
\end{equation*}
This proves the proposition.
\end{proof}

\noindent We now consider the complementary case, namely, points $\tau\in\mathbb{H}$ such that $j(\tau)\in\overline{\mathbb{Q}}$. 
\begin{proof}[Proof of Theorem \ref{general-mero-values-thm}(b)]
Fix a $\tau\in\mathbb{H}$ such that $j(\tau)\in\overline{\mathbb{Q}}$. From Lemma \ref{omegatau}, we get a transcendental number $\omega_\tau$ such that $g_2(\omega_\tau L_\tau)$ and $g_3(\omega_\tau L_\tau)$ are both algebraic. Moreover, we have the formulae
$$
\begin{aligned}
E_4(\tau)&=\frac{3}{4\pi^4}\, g_2(L_\tau)=\frac{3}{4}\left(\frac{\omega_\tau}{\pi} \right)^4 g_2(\omega_\tau L_\tau),\\ E_6(\tau)&=\frac{27}{8\pi^6}\,g_3(L_\tau)=\frac{27}{8}\left(\frac{\omega_\tau}{\pi}\right)^6g_3(\omega_\tau L_\tau).
\end{aligned}
$$
From Nesterenko's Theorem \ref{nest}, we know that at most one of $E_4(\tau)$ and $E_6(\tau)$ is algebraic. The above formulae imply that if ${\omega_{\tau}}/{\pi}\in\overline{\mathbb{Q}}$, then both $E_4(\tau)$ and $E_6(\tau)$ are algebraic, which is a contradiction. Hence, ${\omega_\tau}/{\pi}$ is a transcendental number. Besides, we have 
$$
\begin{aligned}
\Delta(\tau)=\frac{E_4(\tau)^3-E_6(\tau)^2}{1728}&=\frac{1}{4^6}\left(\frac{\omega_\tau}{\pi}\right)^{12}\bigg( g_2(\omega_\tau L_\tau)^3-27 g_3(\omega_\tau L_\tau)^2 \bigg)\\
&=\frac{1}{4^6}\left(\frac{\omega_\tau}{\pi}\right)^{12}\Delta_0(\omega_\tau L_\tau).  
\end{aligned}
$$
As $g_2(\omega_\tau L_\tau)$ and $g_3(\omega_\tau L_\tau)$ both are algebraic, the number $\Delta_0(\omega_\tau L_\tau)\in\overline{\mathbb{Q}}\setminus \{0\}$. Since ${\omega_\tau}/{\pi}$ is transcendental and $\Delta(\tau)$ is a non-zero algebraic multiple of $\left({\omega_\tau}/{\pi}\right)^{12}$, we deduce that $\Delta(\tau)$ is transcendental. \\

Consider the modular function $$g(\tau):=\frac{f(\tau)^{12}}{\Delta(\tau)^k}\quad \forall \tau\in\mathbb{H},$$ which lies in $ \mathcal{F}_{N,\overline{\mathbb{Q}}}$. Since $j(\tau)\in\overline{\mathbb{Q}}$, Theorem \ref{algemf} implies that $g(\tau)$ is algebraic, say $\alpha$. Thus, we get that
\begin{equation}\label{forms-value-j-alg}
f(\tau)^{12}= \alpha \, \cdot \, \Delta(\tau)^k=\frac{\alpha }{4^6} \, \cdot \, \Delta_0(\omega_\tau L_\tau)^k \, \left(\frac{\omega_\tau}{\pi} \right)^{12 k}. 
\end{equation}
This shows that if $f(\tau) \neq 0$, 
 then it is a non-zero algebraic multiple of $\left({\omega_\tau}/{\pi} \right)^{k}$ and hence, is transcendental.
\end{proof}
\bigskip

\subsection{Values of quasi-modular forms}
Let $\widetilde{M}_{k,\overline{\mathbb{Q}}}^{(p)}(\Gamma)$ denote the set of all meromorphic quasi-modular forms of weight $k$ and depth $p(>0)$ for $\Gamma$ with algebraic Fourier coefficients. We study their values at the points $\tau\in\mathbb{H}$, where exactly one of $e^{2\pi i\tau}$ and $j(\tau)$ is algebraic.\\

\begin{proof}[Proof of Theorem \ref{quasi-mdf-thm}(a)]
By Theorem \ref{quasistructure}, we can write $\widetilde{f}$ in the form $$\widetilde{f}=\sum_{r=0}^p f_r E_2^r,\quad f_r\in M_{k,N,\overline{\mathbb{Q}}}.$$ 
Suppose that both $e^{2\pi i\tau}, \widetilde f(\tau)\in \overline{\mathbb{Q}}$. By Corollary \ref{algi-cor}, we know that each number $f_r(\tau)$ is algebraically dependent with $E_4(\tau), E_6(\tau)$. Thus, we get that $E_2(\tau)$ is algebraic over $\overline{\mathbb{Q}}(E_4(\tau), E_6(\tau))$. This implies that
\begin{equation*}
    \text{trdeg}_{\,\mathbb{Q}} \, \,\overline{\mathbb{Q}} \bigg( e^{2 \pi i \tau}, E_2(\tau), E_4(\tau), E_6(\tau) \bigg) \leq 2,
\end{equation*}
contradicting Nesterenko's theorem \ref{nest}. This proves the claim.
\end{proof}

To study values of quasi-modular forms in the complementary case (i.e., for $j(\tau)\in\overline{\mathbb{Q}}$), we need the following lemma.
\begin{lemma}\label{peri}
 Let $\tau\in \mathbb{H}$ be such that $j(\tau)$ is algebraic and $\omega_{\tau}$ be the transcendental number determined in Lemma \ref{omegatau} such that $g_2(\omega_{\tau}L_{\tau})$ and $g_3(\omega_{\tau}L_{\tau)}$ are algebraic. Let $\eta_2(\omega_{\tau}) := \eta_2(\omega_{\tau} L_{\tau})$ be the quasi-period 
 \begin{equation*}
     \eta_2( \omega_{\tau} L_{\tau}) = \zeta_{\omega_{\tau}L_{\tau}}(z+\omega_{\tau}) - \zeta_{\omega_{\tau}L_{\tau}}(z).
 \end{equation*}
 Then $\frac{\omega_\tau}{\pi}$ and $\frac{\eta_2}{\pi}$ are algebraically independent over $\overline{\mathbb{Q}}$. 
\end{lemma}
\begin{proof}
From the definition of the Weierstrass zeta-function, one gets
 $   
{\zeta}_{\omega_\tau L_\tau}(z)= \frac{1}{\omega_\tau}\cdot \zeta_{L_\tau}\left(\frac{z}{\omega_\tau}\right).
$ 
Using the definition of a quasi-period and the identity \eqref{quasi-period}, we obtain
$$
\eta_2= \zeta_{\omega_\tau L_\tau}(z+\omega_\tau)-\zeta_{\omega_\tau L_\tau}(z)= \frac{1}{\omega_\tau}\,\eta_2(L_\tau)=\frac{1}{\omega_\tau}\,G_2(\tau)= \frac{1}{3}\cdot\frac{\pi^2}{ \omega_\tau}\,E_2(\tau).
$$
Thus, we get the following formulae 
$$
E_2(\tau)=3\,\frac{\omega_\tau}{\pi}\frac{\eta_2}{\pi}, \quad E_4(\tau)=\frac{3}{4}\left(\frac{\omega_\tau}{\pi} \right)^4 g_2(\omega_\tau L_\tau),\quad E_6(\tau)=\frac{27}{8}\left(\frac{\omega_\tau}{\pi}\right)^6g_3(\omega_\tau L_\tau).$$
The above formulae imply that $E_2(\tau)$, $E_4(\tau)$ and $E_6(\tau)$ are algebraic over $\overline{\mathbb{Q}}\left(\frac{\omega_\tau}{\pi}, \frac{\eta_2}{\pi}\right)$. Suppose that $\frac{\omega_\tau}{\pi}, \frac{\eta_2}{\pi}$ are algebraically dependent. This implies that $$ \text{trdeg}_{\,\mathbb{Q}} \, \,\overline{\mathbb{Q}} \bigg(  E_2(\tau), E_4(\tau), E_6(\tau) \bigg)=1,$$ which contradicts Nesterenko's Theorem \ref{nest}.
\end{proof}

Lemma \ref{peri} and Theorem \ref{quasistructure} together allow us to describe values of quasi-modular forms explicitly at the points $\tau\in\mathbb{H}$ where the $j$-function is algebraic.
\begin{proof}[Proof of Theorem \ref{quasi-mdf-thm}(b)]
From Theorem \ref{quasistructure}, we have the expression $$\widetilde{f}(\tau)=\sum_{r=0}^p f_r(\tau) E_2(\tau)^r,\quad \text{where } f_r\in M^m_{k,N,\overline{\mathbb{Q}}}.$$
Writing the value $f_r(\tau)$ as in \eqref{forms-value-j-alg} for each coefficient $f_r$, and using the above formula for $E_2$, we obtain
\begin{equation}\label{last}
\widetilde{f}(\tau)=\sum_{r=0}^p c_r\left(\frac{\omega_\tau}{\pi}\right)^{k-2 r}\left(\frac{\omega_\tau}{\pi}\cdot\frac{\eta_2}{\pi}\right)^r=\sum_{r=0}^p\, c_r\, \left(\frac{\omega_\tau}{\pi}\right)^{k-r}\, \left(\frac{\eta_2}{\pi}\right)^r,
 \end{equation}
where each $c_r$ is an algebraic number. If  $c_r \neq 0$ for some $r$ satisfying $0 \leq r \leq p$, then \eqref{last} gives a non-trivial algebraic relation among $\frac{\omega_\tau}{\pi}$ and $\frac{\eta_2}{\pi}$. This contradicts Lemma \ref{peri}. Therefore, the number $\widetilde{f}(\tau)$ is either zero, precisely when each $\widetilde{f}(\tau)\in\overline{\mathbb{Q}}\setminus \{0\}$, or is transcendental.     
\end{proof}

\bigskip

\subsection{Algebraic independence of special values}

\begin{proof}[Proof of Theorem \ref{tr-deg-thm}]
From Theorem \ref{mdft..}, we know that for any $\tau \in \mathbb{H}$ which is not a pole of $g$, $g(\tau)$ is algebraic over $\overline{\mathbb{Q}}(j(\tau))$. Moreover, Corollary \ref{algi-cor} gives that $f(\tau)$ is algebraic over $\overline{\mathbb{Q}}(E_4(\tau), E_6(\tau))$. Since $\widetilde{f}(\tau)=\sum_{r=0}^p f_r(\tau) E_2(\tau)^r$, we have that $\widetilde{f}(\tau)$ is algebraic over $${\overline{\mathbb{Q}}(E_2(\tau), E_4(\tau), E_6(\tau)).}$$ Hence, we get that
\begin{align*}
    \text{trdeg}_{\, \mathbb{Q}} \, \, \overline{\mathbb{Q}} \bigg(e^{2 \pi i \tau}, g(\tau), f(\tau), \widetilde{f}(\tau) \bigg) & = \text{trdeg}_{\, \mathbb{Q}} \, \, \overline{\mathbb{Q}} \bigg(e^{2 \pi i \tau}, j(\tau), E_2(\tau), E_4(\tau), E_6(\tau) \bigg) \\
    & = \text{trdeg}_{\, \mathbb{Q}} \, \, \overline{\mathbb{Q}} \bigg(e^{2 \pi i \tau}, E_2(\tau), E_4(\tau), E_6(\tau) \bigg) \geq 3,
\end{align*}
by Nesterenko's theorem \ref{nest}. This establishes the claim.
\end{proof}

\section*{Acknowledgment}
The authors thank Prof. Michel Waldschmidt and Dr. Veekesh Kumar for insightful suggestions on an earlier version of the paper. They are also grateful to the referee for helpful comments.

\end{document}